\documentclass[a4paper,12pt]{amsart}
\usepackage{amsfonts}
\usepackage{amssymb}
\usepackage{ifthen}
\usepackage{amscd}
\usepackage{amsxtra}
\usepackage{graphicx}
\usepackage{color}
\nonstopmode \numberwithin{equation}{section}
\newtheorem{thm}{Theorem}
\newtheorem{lem}{Lemma}
\newtheorem{cor}{Corollary}
\newtheorem{prop}{Proposition}

\newtheorem{cl}{Claim}
\newtheorem{ca}{Case}
\newtheorem{sca}{Subcase}
\newtheorem{scl}{Subclaim}
\newtheorem{conj}{Conjecture}

\theoremstyle{definition}
\newtheorem{defn}{Definition}

\newtheorem{op}[equation]{Open Problem}
\newtheorem{ques}[equation]{Question}
\newtheorem{rem}{Remark}[section]
\newtheorem{exam}[equation]{Example}

\newcounter {own}
\def\theown {\thesection     .\arabic{own}}

\newenvironment{pf}[1][]{%
 \vskip 3mm
 \noindent
 \ifthenelse{\equal{#1}{}}%
 {{\slshape Proof. }}%
 {{\slshape #1.} }%
 }%
{\qed\bigskip}

\newcounter{alphabet}
\newcounter{tmp}
\newenvironment{Thm}[1][]{\refstepcounter{alphabet}%
\bigskip%
\noindent%
{\bf Theorem \Alph{alphabet}}%
\ifthenelse{\equal{#1}{}}{}{ (#1)}%
{\bf .} \itshape}{\vskip 8pt}

\makeatletter
\newcommand{\Ref}[1]{\@ifundefined{r@#1}{}{\setcounter{tmp}{\ref{#1}}\Alph{tmp}}}
\makeatother

\newenvironment{Lem}[1][]{\refstepcounter{alphabet}%
\bigskip%
\noindent%
{\bf Lemma \Alph{alphabet}}%
{\bf .} \itshape}{\vskip 8pt}

\newcommand{\IR}{{\mathbb R}}

\newcommand{\Aut}{{\operatorname{Aut}}}

\def\be{\begin{equation}}
\def\ee{\end{equation}}

\newcommand{\bee}{\begin{enumerate}}
\newcommand{\eee}{\end{enumerate}}

\newcommand{\blem}{\begin{lem}}
\newcommand{\elem}{\end{lem}}
\newcommand{\bthm}{\begin{thm}}
\newcommand{\ethm}{\end{thm}}
\newcommand{\bcor}{\begin{cor}}
\newcommand{\ecor}{\end{cor}}
\newcommand{\beg}{\begin{exam}}
\newcommand{\eeg}{\end{exam}}
\newcommand{\begs}{\begin{examples}}
\newcommand{\eegs}{\end{examples}}
\newcommand{\bdefe}{\begin{defn}}
\newcommand{\edefe}{\end{defn}}
\newcommand{\bprob}{\begin{prob}}
\newcommand{\eprob}{\end{prob}}
\newcommand{\bques}{\begin{ques}}
\newcommand{\eques}{\end{ques}}
\newcommand{\bei}{\begin{itemize}}
\newcommand{\eei}{\end{itemize}}
\newcommand{\bcon}{\begin{conj}}
\newcommand{\econ}{\end{conj}}
\newcommand{\bop}{\begin{op}}
\newcommand{\eop}{\end{op}}

\newcommand{\bca}{\begin{ca}}
\newcommand{\eca}{\end{ca}}
\newcommand{\bsca}{\begin{sca}}
\newcommand{\esca}{\end{sca}}

\newcommand{\bcl}{\begin{cl}}
\newcommand{\ecl}{\end{cl}}

\newcommand{\bscl}{\begin{scl}}
\newcommand{\escl}{\end{scl}}

\newcommand{\bcons}{\begin{conjs}}
\newcommand{\econs}{\end{conjs}}
\newcommand{\bprop}{\begin{propo}}
\newcommand{\eprop}{\end{propo}}
\newcommand{\br}{\begin{rem}}
\newcommand{\er}{\end{rem}}
\newcommand{\brs}{\begin{rems}}
\newcommand{\ers}{\end{rems}}
\newcommand{\bo}{\begin{obser}}
\newcommand{\eo}{\end{obser}}
\newcommand{\bos}{\begin{obsers}}
\newcommand{\eos}{\end{obsers}}
\newcommand{\bpf}{\begin{pf}}
\newcommand{\epf}{\end{pf}}
\newcommand{\ba}{\begin{array}}
\newcommand{\ea}{\end{array}}
\newcommand{\beq}{\begin{eqnarray}}
\newcommand{\beqq}{\begin{eqnarray*}}
\newcommand{\eeq}{\end{eqnarray}}
\newcommand{\eeqq}{\end{eqnarray*}}

\newcommand{\ds}{\displaystyle}

\newcounter{minutes}
\divide\time by 60
\newcounter{hours}
\multiply\time by 60 \addtocounter{minutes}{-\time}

\begin{document}
\bibliographystyle{amsplain}
\title [] {Schwarz's Lemmas for mappings satisfying
Poisson's equation }

\def\thefootnote{}
\footnotetext{ \texttt{\tiny File:~\jobname .tex,
          printed: \number\day-\number\month-\number\year,
          \thehours.\ifnum\theminutes<10{0}\fi\theminutes}
} \makeatletter\def\thefootnote{\@arabic\c@footnote}\makeatother

\author{Shaolin Chen}
\address{Sh. Chen,    College of Mathematics and
Statistics, Hengyang Normal University, Hengyang, Hunan 421008,
People's Republic of China} \email{mathechen@126.com}

\author{Saminathan Ponnusamy 
}
\address{S. Ponnusamy, Stat-Math Unit,
Indian Statistical Institute (ISI), Chennai Centre,
110, Nelson Manickam Road,
Aminjikarai, Chennai, 600 029, India. }
\email{samy@isichennai.res.in, samy@iitm.ac.in}


\begin{abstract}
For $n\geq3$, $m\geq1$ and a given continuous function
$g:~\Omega\rightarrow\mathbb{R}^{m}$, we establish some Schwarz type
lemmas for mappings $f$ of $\Omega$ into $\mathbb{R}^{m}$ satisfying
the {\rm PDE}: $\Delta f=g$, where $\Omega$ is a subset of
$\mathbb{R}^{n}$. Then we apply these
 results to obtain a Landau type theorem.
\end{abstract}



\maketitle \pagestyle{myheadings}
\markboth{Sh. Chen and S. Ponnusamy} {Schwarz's Lemmas for mappings satisfying Poisson's
equation}

\section{Preliminaries and statements of main results}\label{csw-sec1}
\subsection{Notations}
For $n\geq2$, let $\mathbb{R}^{n}
$ denote the usual real vector space of dimension $n$. Sometimes it is convenient to
identify each point $x=(x_{1},\ldots,x_{n})\in\mathbb{R}^{n}$ with an $n\times 1$ column matrix
so that
$$x=\left(\begin{array}{cccc}
x_{1} \\
\vdots \\
 x_{n}
\end{array}\right).
$$
For $w=(w_{1},\ldots,w_{n})$ and $x\in\mathbb{R}^{n}$, we define the
Euclidean inner product $\langle \cdot ,\cdot \rangle$ by
$\langle x,w\rangle=x_{1}w_{1}+\cdots+x_{n}w_{n}
$
so that the Euclidean length of $x$ is defined by
$$|x|=\langle x,x\rangle^{1/2}=(|x_{1}|^{2}+\cdots+|x_{n}|^{2})^{1/2}.
$$
Denote a ball in $\mathbb{R}^{n}$ with center $w\in\mathbb{R}^{n}$
and radius $r$ by
$$\mathbb{B}^{n}(w,r)=\{x\in\mathbb{R}^{n}:\, |x-w|<r\}.
$$
In particular, $\mathbb{B}^{n}$ and $\mathbb{S}^{n-1}$ denote the unit ball
$\mathbb{B}^{n}(0,1)$  and the unit sphere in $\mathbb{R}^{n}$, respectively.
Set $\mathbb{D}=\mathbb{B}^2$, the open unit
disk in the complex plane $\mathbb{C}\cong \mathbb{R}^{2}$. For
$k\in\mathbb{N}_0:=\mathbb{N}\cup\{0\}$ and $m\in\mathbb{N}_0$, we
denote by $\mathcal{C}^{k}(\Omega_{1},\Omega_{2})$ the set of all
$k$-times continuously differentiable functions from $\Omega_{1}$
into $\Omega_{2}$, where $\Omega_{1}$ and $\Omega_{2}$ are subsets
of $\mathbb{R}^{n}$ and $\mathbb{R}^{m}$, respectively. In
particular, let
$\mathcal{C}(\Omega_{1},\Omega_{2}):=\mathcal{C}^{0}(\Omega_{1},\Omega_{2})$,
the set of all continuous functions of $\Omega_{1}$ into
$\Omega_{2}$. For
$f=(f_{1},\ldots,f_{m})\in\mathcal{C}^{1}(\Omega_{1},\Omega_{2})$,
we denote the derivative $D_{f}$ of $f$ by
$$D_{f}=\left(\begin{array}{cccc}
\ds D_{1}f_{1}\; \cdots\;
 D_{n}f_{1}\\[4mm]
\vdots\;\; \;\;\cdots\;\;\;\;\vdots \\[2mm]
 \ds
D_{1}f_{m}\; \cdots\;
 D_{n}f_{m}
\end{array}\right), \quad D_{j}f_{i}(x)=\frac{\partial f_{i}(x)}{\partial x_j}.
$$
In particular, if $n=m$, the Jacobian of $f$ is defined by
$J_{f}=\det D_{f}$ and the Laplacian of
$f\in\mathcal{C}^{1}(\Omega_{1},\Omega_{2})$ is defined by
$$\Delta f=\sum_{k=1}^{n}D_{kk}f.
$$
For an $m\times n$ matrix $A$, the operator norm of $A$ is defined by
$$|A|=\sup_{x\neq 0}\frac{|Ax|}{|x|}=\max\{|A\theta|:\,
\theta\in\mathbb{S}^{n-1}\}
$$
where $m\geq1$ and $n\geq2$.

\subsection{Poisson equation and Schwarz lemma}
For $x,y\in\mathbb{R}^{n}\backslash\{0\}$, we define $x'=x/|x|$, $y'=y/|y|$ and let
$$[x,y]:=\left|y|x|-x'\right|=\left|x|y|-y'\right|.
$$
Also, for $x,y\in\mathbb{B}^{n}$ with $x\neq y$ and $|x|+|y|\neq0$,
we use $G(x,y)$ to denote the {\it Green function}:
\be\label{eq-ex0}
G(x,y)=c_{n}\left(\frac{1}{|x-y|^{n-2}}-\frac{1}{[x,y]^{n-2}}\right),
\ee where $c_{n}=1/[(n-2)\omega_{n-1}]$ and
$\omega_{n-1}=2\pi^{\frac{n}{2}}/\Gamma\big(\frac{n}{2}\big)$
denotes the {\it Hausdorff measure} of $\mathbb{S}^{n-1}$. The {\it
Poisson kernel} $P:\,\mathbb{B}^{n}\times
\mathbb{S}^{n-1}\rightarrow \IR$ is defined by
$$P(x,\zeta)=\frac{1-|x|^{2}}{|x-\zeta|^{n}}.
$$
We write
$$\nabla =\left ( \frac{\partial }{\partial x_1}, \ldots, \frac{\partial }{\partial x_n} \right )
$$
and for a vector valued function $f=(f_1, \ldots, f_m)$, we define
the directional derivative $\frac{\partial f}{\partial
\overrightarrow{n}}$ by a componentwise interpretation:
$$ \frac{\partial f}{\partial \overrightarrow{n}}(x)=\left ( \langle \nabla f_1(x), x'\rangle, \ldots, \ldots, \langle \nabla f_m(x), x'\rangle \right ).
$$
For a given bounded integrable function $\psi:~\mathbb{S}^{n-1}\rightarrow\mathbb{R}^{m}$ and
$g\in\mathcal{C}(\mathbb{B}^{n},\mathbb{R}^{m})$, the solution of the {\it Poisson equation}
\be\label{eq-1}
\begin{cases} \displaystyle \Delta f=g
& \mbox{ in } \mathbb{B}^{n}\\
\displaystyle f=\psi &\mbox{ in }\, \mathbb{S}^{n-1}
\end{cases}
\ee is given by \be\label{eq-p}
f(x)=\mathcal{P}_{\psi}(x)-\mathcal{G}_{g}(x), \ee where
\be\label{eq-p1a}
\mathcal{P}_{\psi}(x)=\int_{\mathbb{S}^{n-1}}P(x,\zeta)\psi(\zeta)d\sigma(\zeta)
~\mbox{ and }~\mathcal{G}_{g}(x)=
\int_{\mathbb{B}^{n}}G(x,y)g(y)dV(y) \ee for $x\in\mathbb{B}^{n}$.
Here $d\sigma$ denotes the normalized Lebesgue surface measure on
$\mathbb{S}^{n-1}$ and $dV$ is the Lebesgue volume measure on
$\mathbb{B}^{n}$. It is well known that if $\psi$ and $g$ are
continuous in $\mathbb{S}^{n-1}$ and in $\overline{\mathbb{B}^{n}}$,
respectively, then $f=\mathcal{P}_{\psi}-\mathcal{G}_{g}$ has a
continuous extension $\tilde{f}$ to the boundary, and
$\tilde{f}=\psi$ in $\mathbb{S}^{n-1}$ (see \cite[p.~118--120]{Ho}
or \cite{K4,K3,K1}).

The classical Schwarz lemma states that a holomorphic function $f$
from $\mathbb{D}$ into itself with $f(0)=0$ satisfies
$|f(z)|\leq|z|$ for all $z\in\mathbb{D}$. It is well-known that the
Schwarz lemma has become a crucial theme in a lot of branches of
mathematical research for more than a hundred years to date. For
$n\geq3$,  the classical Schwarz lemma of harmonic
mappings in $\mathbb{B}^{n}$ infers that if $f$ is a
harmonic mapping of $\mathbb{B}^{n}$ into itself satisfying
$f(0)=0,$ then
$$|f(x)|\leq U(rN),
$$
where $r=|x|$, $N=(0,\ldots,0,1)$ and $U$ is
a harmonic function of $\mathbb{B}^{n}$ into $[-1,1]$ defined by
$$U(x)=\mathcal{P}_{(\mathcal{X}_{S^{+}}-\mathcal{X}_{S^{-}})}(x).
$$
Here $\mathcal{X}$ is the indicator function,
$S^{+}=\{x=(x_{1},\ldots,x_{n})\in\mathbb{S}^{n-1}:~x_{n}\geq0\}$
and
$S^{-}=\{x=(x_{1},\ldots,x_{n})\in\mathbb{S}^{n-1}:~x_{n}\leq0\}$
(see \cite{ABR}). For the case $n=2,$ we refer to \cite{CK,CV,He,
K0}.
In \cite{K5}, Kalaj proved the following result for harmonic mappings $f$ of $\mathbb{B}^{n}$
into itself:
\be\label{eq-K}
\left|f(x)-\frac{1-|x|^{2}}{(1+|x|^{2})^{\frac{n}{2}}}f(0)\right|\leq U(|x|N).
\ee
\subsection{Main results and Remarks}
The first aim of the paper is to extend the result \eqref{eq-K} to mappings satisfying the Poisson equation.
More precisely, we shall prove the following.

\begin{thm}\label{thm-1}
Let $n\geq3$, $m\geq1$ and $g\in\mathcal{C}(\overline{\mathbb{B}^{n}},\mathbb{R}^{m})$. If
$f\in\mathcal{C}^{2}(\mathbb{B}^{n},\mathbb{R}^{m})\cap\mathcal{C}(\mathbb{S}^{n-1},\mathbb{R}^{m})$
satisfies $\Delta f=g$, then, for $x\in\overline{\mathbb{B}^{n}}$,
\be\label{eq-thm1}
\left|f(x)-\frac{1-|x|^{2}}{(1+|x|^{2})^{\frac{n}{2}}}\mathcal{P}_{f}(0)\right|
\leq\|\mathcal{P}_{f}\|_{\infty}U(|x|N)
+\frac{\|g\|_{\infty}}{2n}(1-|x|^{2}),\ee where
$$\mathcal{P}_{f}(x)=\int_{\mathbb{S}^{n-1}}P(x,\zeta)f(\zeta)d\sigma(\zeta),~\|f\|_{\infty}=
\sup_{x\in\mathbb{B}^{n}}|f(x)|~\mbox{and}~\|g\|_{\infty}=\sup_{x\in\mathbb{B}^{n}}|g(x)|.
$$
If we choose $g(x)=(-2nM,0,\ldots,0)$ and
$f(x)=(M(1-|x|^{2}),0,\ldots,0)$ for
$x\in\overline{\mathbb{B}^{n}}$, then
 the inequality
{\rm(\ref{eq-thm1})} is sharp in $\mathbb{S}^{n-1}\cup\{0\}$, where
$M>0$ is a constant.
\end{thm}

By Theorem \ref{thm-1}, we give an explicit estimate as follows.

\begin{cor}\label{cor-1}
Let $n\geq3$, $m\geq1$ and $g\in\mathcal{C}(\overline{\mathbb{B}^{n}},\mathbb{R}^{m})$. If
$f\in\mathcal{C}^{2}(\mathbb{B}^{n},\mathbb{R}^{m})\cap\mathcal{C}(\mathbb{S}^{n-1},\mathbb{R}^{m})$
satisfies $\Delta f=g$, then, for $x\in\overline{\mathbb{B}^{n}}$,
\be\label{eq-thm1a}
\left|f(x)-\frac{1-|x|}{(1+|x|)^{n-1}}\mathcal{P}_{f}(0)\right|
\leq\|\mathcal{P}_{f}\|_{\infty}\left[1-\frac{1-|x|}{(1+|x|)^{n-1}}\right]
+\frac{\|g\|_{\infty}}{2n}(1-|x|^{2}).
\ee
If we choose $g(x)=(-2nM,0,\ldots,0)$ and $f(x)=(M(1-|x|^{2}),0,\ldots,0)$ for
$x\in\overline{\mathbb{B}^{n}}$, then  the inequality {\rm(\ref{eq-thm1a})} is sharp
in $\mathbb{S}^{n-1}\cup\{0\}$, where $M>0$ is a constant.
\end{cor}

There is a classical Schwarz lemma at the boundary, which is as
follows and the same may be obtained from standard texts. See for instance,
\cite{G} or \cite[p.~249, Corollary~6.62]{pon}.

\begin{Thm}
\label{Thm-B}
Let $f$ be a holomorphic function from
$\mathbb{D}$ into itself. If $f$ is holomorphic at $z=1$ with
$f(0)=0$ and $f(1)=1$, then $f'(1)\geq1$. Moreover, the inequality
is sharp.
\end{Thm}

Theorem \Ref{Thm-B} has been generalized in various forms. For
example, Krantz \cite{Kra} explored many versions of the Schwarz
lemma at the boundary point of a domain in $\mathbb{C}$, and
reviewed several results of many authors. A natural response to
these results is to ask whether we can extend the Schwarz lemma at
the boundary to higher dimensional cases. This form of the Schwarz lemma has attracted much
attention, see \cite{BK,LWT,LT} for holomorphic functions, and see
\cite{ABR,Bu, K5,K6,MM} for harmonic functions. In the following, by
using Theorem \ref{thm-1}, we establish a Schwarz lemma at the
boundary for mappings satisfying Poisson's equation, which is also a
generalization of Theorem \Ref{Thm-B}.

\begin{thm}\label{thm-2} For $n\geq3$, $m\geq1$ and a given
$g\in\mathcal{C}(\overline{\mathbb{B}^{n}},\mathbb{R}^{m})$, let
$f\in\mathcal{C}^{2}(\mathbb{B}^{n},\mathbb{R}^{m})\cap\mathcal{C}(\mathbb{S}^{n-1},\mathbb{R}^{m})$
be a mapping of $\mathbb{B}^{n}$ into itself satisfying $\Delta f=g,$
where
$$\|g\|_{\infty}<\frac{nA_{n}}{\big(1+\frac{1}{2^{\frac{n}{2}}}\big)} ~\mbox{ and }~ A_{n}=\frac{n!\big[1+n-(n-2)F\big(\frac{1}{2},1;\frac{n+3}{2};-1\big)\big]}{2^{\frac{3n}{2}}
\Gamma\big(\frac{n+1}{2}\big)\Gamma\big(\frac{n+3}{2}\big)}.
$$
If $f(0)=0$ and $\lim_{r\rightarrow1^{-}}|f(r\zeta)|=1$ for some $\zeta\in\mathbb{S}^{n-1}$, then
\be\label{eq-Sch}
\liminf_{r\rightarrow1^{-}}\frac{|f(\zeta)-f(r\zeta)|}{1-r}\geq
A_{n}-\frac{\|g\|_{\infty}}{n}\left(1+\frac{1}{2^{\frac{n}{2}}}\right).\ee
In particular, if $\|g\|_{\infty}=0$, then the estimate of
{\rm(\ref{eq-Sch})} is sharp.
\end{thm}

The definition of the classical hypergeometric function $F(a,b;c;z)$ is given in Section \ref{csw-sec2}.
Next, we state a Schwarz-Pick type lemma, which is a generalization
of \cite[Theorem 1.3]{K6}, \cite[Theorem 2.12]{MM} and
\cite[Corollary 2.2]{K5}.

\begin{thm}\label{thm-3}
Let $n\geq3$ and $m\geq1$. If
$f\in\mathcal{C}^{2}(\mathbb{B}^{n},\mathbb{R}^{m})\cap\mathcal{C}(\mathbb{S}^{n-1},\mathbb{R}^{m})$
satisfies $\Delta f=g$ for a given $g\in\mathcal{C}(\overline{\mathbb{B}^{n}},\mathbb{R}^{m})$, then
$$|D_{f}(x)|\leq\frac{\|\mathcal{P}_{f}\|_{\infty}}{1-r}\sup_{\gamma>0}C(\gamma,r)+
\frac{n}{n+1}\|g\|_{\infty} ~\mbox{ for $x\in\mathbb{B}^{n}$},
$$
where $r=|x|$ and
$$C(\gamma,r)=\frac{4\omega_{n-3}}{\omega_{n-1}}
\frac{2^{n-1}}{(1+r)^{n-1}}\frac{1}{\sqrt{1+\gamma^{2}}}\int_{0}^{1}\frac{\Psi_{r}(\gamma
t)+\Psi_{r}(-\gamma t)}{\sqrt{(1-t^{2})^{4-n}}}dt<\infty.
$$
Here
$$\Psi_{r}(\gamma)=\int_{0}^{\frac{\gamma+\sqrt{\gamma^{2}+1-\alpha^{2}(r)}}{1-\alpha(r)}}
\frac{n-\beta(r)+n\gamma
\rho-\beta(r)\rho^{2}}{(1+\rho^{2})^{\frac{n}{2}+1}(1+\tau^{2}(r)\rho^{2})^{\frac{n}{2}-1}}\rho^{n-2}d\rho,
$$
$$\tau(r)=\frac{1-r}{1+r},~\alpha(r)=\frac{r(n-2)}{n}~\mbox{and}~\beta(r)=\frac{[n-(n-2)r]}{2}.
$$
\end{thm}

\br
We observe that if $n=4$ in Theorem \ref{thm-3}, then
$\sup_{\gamma>0}C(\gamma,r)$ can be replaced by
$$\frac{\left[r\sqrt{4-r^{2}}(2+r^{2})+4(1-r^{2})\arctan\big(r\frac{\sqrt{4-r^{2}}}{r^{2}-2}\big)\right]}{\pi(1+r)r^{3}},
$$
where $r=|x|$ (see \cite[Theorem 1.3]{K6}). In particular, using
(\ref{eq-K}), we can obtain an explicit estimate for
$|D_{\mathcal{P}_{f}}(x)|$. But it is not going to yield a better
bound than the bound stated in Theorem \ref{thm-3}. Indeed, for any fixed
$x\in\mathbb{B}^{n}$, consider the function
$$\nu(y)=\mathcal{P}_{f}\big(x+(1-|x|)y\big).
$$
Now, if we apply (\ref{eq-K}) to $\nu/\|\mathcal{P}_{f}\|_{\infty}$,
then we deduce that
$$\left|\frac{\mathcal{P}_{f}\big(x+(1-|x|)y\big)-\mathcal{P}_{f}(x)}{|y|}-
\frac{\left[(1-|y|^{2})(1+|y|^{2})^{-n/2}-1\right]}{|y|}\mathcal{P}_{f}(x)\right|\leq\|\mathcal{P}_{f}\|_{\infty}\frac{U(|y|N)}{|y|},
$$
which, together with the fact that
$$\lim_{t\rightarrow0^{+}}\frac{(1-t^{2})(1+t^{2})^{-n/2}-1}{t}=0,
$$
implies \be\label{eq18c}
(1-|x|)|D_{\mathcal{P}_{f}}(x)|\leq\|\mathcal{P}_{f}\|_{\infty}\frac{\partial
U(rN)}{\partial
r}\Big|_{r=0}=\frac{2\omega_{n-1}}{V(\mathbb{B}^{n})}\|\mathcal{P}_{f}\|_{\infty}=2n\|\mathcal{P}_{f}\|_{\infty},\ee
where $V(\mathbb{B}^{n})$ is the volume of $\mathbb{B}^{n}$.
Finally, by  (\ref{eq18c}) and (\ref{eq-15c}),   we conclude that
\be\label{eq-19cp}|D_{f}(x)|\leq|D_{\mathcal{P}_{f}}(x)|+|D_{\mathcal{G}_{g}}(x)|\leq\frac{2n\|\mathcal{P}_{f}\|_{\infty}}{1-|x|}+
\frac{n}{n+1}\|g\|_{\infty} ~\mbox{ for $x\in\mathbb{B}^{n}$} \ee as
desired. \hfill $\Box$ \er


There exists a number of articles dealing with Landau type theorems
in geometric function theory and, for general class of functions
without some additional condition(s), there is no Landau's theorem.
See for example \cite{CK,CMPW,CP,CPW-211,CPW-2011,CPW,CPW-2015,CV,W}
and the related references therein. In our next result, we use
Theorems \ref{thm-1} and \ref{thm-3} to establish a Landau type
theorem for mappings satisfying the Poisson equation.

For $n\geq3$ and a given
$g\in\mathcal{C}(\overline{\mathbb{B}^{n}},\mathbb{R}^{n})$, we use
$\mathcal{F}_{g}$ to denote the set of all mappings
$f\in\mathcal{C}^{2}(\mathbb{B}^{n},\mathbb{R}^{n})\cap\mathcal{C}(\mathbb{S}^{n-1},\mathbb{R}^{n})$
satisfying $\Delta f=g$ and $|f(0)|=J_{f}(0)-1=0$. Let
$\mathcal{F}_{g}^{M}$ be the set of all mappings
$f\in\mathcal{F}_{g}$ satisfying $\|f\|_{\infty}+\|g\|_{\infty}\leq
M$, where  $M$ is a positive constant.

\begin{thm}\label{thm-4}
For $n\geq3$ and a given $g\in\mathcal{C}(\overline{\mathbb{B}^{n}},\mathbb{R}^{n})$, let
$f\in\mathcal{F}_{g}^{M}$.  Then there is a positive constant $r_{0}$ depending only on $M$ and $g$ such that
$\mathbb{B}^{n}(0,r_{0})\subset f(\mathbb{B}^{n}).$
\end{thm}

\begin{rem}
We wish to point out that the Landau type theorem fails for the
family $\mathcal{F}_{g}$ without appropriate additional
condition(s). For example, consider $g(x)=(0,\ldots,0, 2n/3)$ and
$f_{k}(x)=(kx_{1},x_{2}/k, x_{3},\ldots,x_{n-1},|x|^{2}/3+x_{n})$
for $k=\{1,2,\ldots\}$, where $n\geq3$ and
$x=(x_{1},\ldots,x_{n})\in\mathbb{B}^{n}$. It is easy to see that
each $f_{k}$ is univalent and $|f_{k}(0)|=J_{f_{k}}(0)-1=0$.
Furthermore, each  $f_{k}(\mathbb{B}^{n})$ contains no ball with
radius bigger than $1/k$. Hence, there does not exist an absolute
constant $r_{0}>0$ which can work for all $k\in\{1,2,\ldots\}$, such
that $\mathbb{B}(0,r_{0}) $ is contained in the range
$f_{k}(\mathbb{B}^{n})$. Although Theorem \ref{thm-4} provides
existence of the Landau-Bloch constant for
$f\in\mathcal{F}_{g}^{M}$, an explicit estimate on the Landau-Bloch
constant is not established.
\end{rem}

The proofs of Theorems \ref{thm-1}, \ref{thm-2}, \ref{thm-3} and Corollary \ref{cor-1} will be
presented in Section \ref{csw-sec2}. Moreover, the proof of Theorem \ref{thm-4} will be
given in Section \ref{csw-sec3}.

\section{The Schwarz Lemmas for mappings satisfying Poisson's equation}\label{csw-sec2}

\subsection{M\"obius Transformations of the Unit Ball}\label{sbcsw-sec2.1}

For $x\in\mathbb{B}^{n}$, the {\it M\"obius transformation} in
$\mathbb{B}^{n}$ is defined by \be\label{eq-ex1}
\phi_{x}(y)=\frac{|x-y|^{2}x-(1-|x|^{2})(y-x)}{[x,y]^{2}},~y\in\mathbb{B}^{n}.
\ee The set of isometries of the hyperbolic unit ball is a {\it
Kleinian subgroup} of all M\"obius transformations of the extended
spaces $\mathbb{R}^{n}\cup\{\infty\}$ onto itself. In the following,
we make use of the {\it automorphism group} $\Aut(\mathbb{B}^{n})$
consisting of all M\"obius transformations of the unit ball
$\mathbb{B}^{n}$ onto itself. We recall the following facts from
\cite{Bea}: For $x\in\mathbb{B}^{n}$ and
$\phi_{x}\in\Aut(\mathbb{B}^{n})$, we have $\phi_{x}(0)=x$,
$\phi_{x}(x)=0$, $\phi_{x}(\phi_{x}(y))=y \in\mathbb{B}^{n}$,
\be\label{II}
|\phi_{x}(y)|=\frac{|x-y|}{[x,y]}, ~1-|\phi_{x}(y)|^{2}=\frac{(1-|x|^{2})(1-|y|^{2})}{[x,y]^{2}}
\ee
and
\be\label{III}
|J_{\phi_{x}}(y)|=\frac{(1-|x|^{2})^{n}}{[x,y]^{2n}}.
\ee

\subsection{Gauss Hypergeometric Functions}\label{sbcsw-sec2.2}

For $a, b, c\in\mathbb{R}$ with $c\neq0, -1, -2, \ldots,$ the {\it
hypergeometric} function is defined by the power series in the variable $x$
$$F(a,b;c;x)=\sum_{k=0}^{\infty}\frac{(a)_{k}(b)_{k}}{(c)_{k}}\frac{x^{k}}{k!},~|x|<1.
$$
Here $(a)_{0}=1$,  $(a)_{k}=a(a+1)\cdots(a+k-1)$ for $k=1, 2, \ldots$, and
generally $(a)_{k}=\Gamma(a+k)/\Gamma(a)$ is the {\it Pochhammer} symbol,
where $\Gamma$ is the {\it Gamma function}. In particular, for $a, b, c>0$ and
$a+b<c$, we have (cf. \cite{PBM}) 
$$F(a,b;c;1)=\lim_{x\rightarrow1}
F(a,b;c;x)=\frac{\Gamma(c)\Gamma(c-a-b)}{\Gamma(c-a)\Gamma(c-b)}<\infty.
$$

The following result is useful in showing one of our main
results of the paper.

\begin{prop}{\rm (\cite{K2}~$\mbox{or}$~\cite[2.5.16(43)]{PBM})}\label{pro-1}
For $\lambda_{1}>1$ and $\lambda_{2}>0$, we have
$$\int_{0}^{\pi}\frac{\sin^{\lambda_{1}-1}t}{(1+r^{2}-2r\cos t)^{\lambda_{2}}}dt=
\mathbf{B}\left(\frac{\lambda_{1}}{2},\frac{1}{2}\right)
F\big(\lambda_{2},\lambda_{2}+\frac{1-\lambda_{1}}{2};\frac{1+\lambda_{1}}{2};r^{2}\big),
$$
where $\mathbf{B}(.,.)$ denotes the beta function and $r\in[0,1)$.
\end{prop}

\subsection{Proofs}
\subsection*{Proof of Theorem \ref{thm-1}}
Let $n\geq3$.
For $x,y\in\mathbb{B}^{n}$ with $x\neq y$ and $|x|+|y|\neq0$, by
\eqref{II}, we have
\beq\label{eq-2}
\left|\frac{1}{|x-y|^{n-2}}-\frac{1}{[x,y]^{n-2}}\right|&=&\frac{1}{|x-y|^{n-2}}\left|1-\frac{|x-y|^{n-2}}{[x,y]^{n-2}}\right|\\
\nonumber&=&\frac{1}{|x-\phi_{x}(z)|^{n-2}}\left(1-|z|^{n-2}\right),
\eeq
where $\phi_{x}\in\Aut(\mathbb{B}^{n})$ and $z=\phi_{x}(y)$. By \eqref{eq-ex1}, direct calculation shows that
\begin{eqnarray*}
x-\phi_{x}(z)=\frac{x[x,z]^{2}-|x-z|^{2}x+(1-|x|^{2})(z-x)}{[x,z]^{2}}
=\frac{(z-x|z|^{2})(1-|x|^{2})}{[x,z]^{2}},
\end{eqnarray*}
which gives \be\label{eq-3}
|x-\phi_{x}(z)|=\frac{|z|(1-|x|^{2})}{[x,z]}. \ee By (\ref{eq-2})
and (\ref{eq-3}), we obtain \be\label{eq-4}
\left|\frac{1}{|x-y|^{n-2}}-\frac{1}{[x,y]^{n-2}}\right|=\frac{[x,z]^{n-2}(1-|z|^{n-2})}{|z|^{n-2}(1-|x|^{2})^{n-2}}.
\ee Now, let
$g\in\mathcal{C}(\overline{\mathbb{B}^{n}},\mathbb{R}^{m})$ be
given. Then, by (\ref{eq-p}) with $f$ in place of $\psi$, we have
$f(x)=\mathcal{P}_{f}(x)-\mathcal{G}_{g}(x), $ where
$\mathcal{P}_{f}$ and $\mathcal{G}_{g}$ are defined by
\eqref{eq-p1a}.
It follows that
\beq\label{eq-5}
\nonumber|\mathcal{G}_{g}(x)|&\leq&\int_{\mathbb{B}^{n}}|G(x,y)g(y)|dV(y)\\
\nonumber
&\leq&c_{n}\|g\|_{\infty}\int_{\mathbb{B}^{n}}\left|\frac{1}{|x-y|^{n-2}}-\frac{1}{[x,y]^{n-2}}\right|dV(y),~\mbox{ by \eqref{eq-ex0}},\\
\nonumber &=&c_{n}\|g\|_{\infty}\int_{\mathbb{B}^{n}}
\frac{[x,z]^{n-2}(1-|z|^{n-2})}{|z|^{n-2}(1-|x|^{2})^{n-2}}\frac{(1-|x|^{2})^{n}}{[x,z]^{2n}}dV(z),~\mbox{ by   \eqref{II} and (\ref{eq-4})},\\
\nonumber
&=&c_{n}\|g\|_{\infty}(1-|x|^{2})^{2}\int_{\mathbb{B}^{n}}\frac{(1-|z|^{n-2})|z|^{n+2}}{|z|^{n-2}\big|x|z|^{2}-z\big|^{2+n}}dV(z)\\
&=&\frac{\|g\|_{\infty}(1-|x|^{2})^{2}}{n-2}\int_{0}^{1}dr\int_{\mathbb{S}^{n-1}}\frac{r(1-r^{n-2})}{|rx-\zeta|^{2+n}}d\sigma(\zeta).
\eeq

Using the polar coordinates and Proposition \ref{pro-1}, we obtain
\beq\label{eq-6}\nonumber
\int_{\mathbb{S}^{n-1}}\frac{d\sigma(\zeta)}{|rx-\zeta|^{2+n}}&=&\frac{1}{\int_{0}^{\pi}\sin^{n-2}tdt}\int_{0}^{\pi}
\frac{\sin^{n-2}tdt}{\left(1+r^{2}|x|^{2}-2r|x|\cos
t\right)^{\frac{n+2}{2}}}dt\\ \nonumber
&=&\frac{\Gamma\big(\frac{n}{2}\big)}{\sqrt{\pi}\Gamma\big(\frac{n-1}{2}\big)}\cdot
\frac{\sqrt{\pi}\Gamma\big(\frac{n-1}{2}\big)}{\Gamma\big(\frac{n}{2}\big)}F\Big(\frac{n+2}{2},2;\frac{n}{2};r^{2}|x|^{2}\Big)\\
&=& F\Big(\frac{n+2}{2},2;\frac{n}{2};r^{2}|x|^{2}\Big).
\eeq
By (\ref{eq-5}) and (\ref{eq-6}), we get
\beq\label{eq-7}\nonumber
|\mathcal{G}_{g}(x)|&\leq&\frac{\|g\|_{\infty}(1-|x|^{2})^{2}}{n-2}
\int_{0}^{1}F\Big(\frac{n+2}{2},2;\frac{n}{2};r^{2}|x|^{2}\Big)
r(1-r^{n-2})dr\\ \nonumber
&=&\frac{\|g\|_{\infty}(1-|x|^{2})^{2}}{n-2}\int_{0}^{1}\sum_{k=0}^{\infty}
\frac{\Gamma\big(\frac{n}{2}\big)\Gamma\big(\frac{n}{2}+k+1\big)
\Gamma(2+k)}{\Gamma\big(\frac{n}{2}+1\big)\Gamma\big(\frac{n}{2}+k\big)\Gamma(2)}
\frac{r^{2k+1}(1-r^{n-2})|x|^{2k}}{k!}\\ \nonumber
&=&\frac{\|g\|_{\infty}(1-|x|^{2})^{2}}{n-2}\int_{0}^{1}
\sum_{k=0}^{\infty}\frac{2}{n}\big(k+\frac{n}{2}\big)(k+1)r^{2k+1}(1-r^{n-2})|x|^{2k}dr\\
\nonumber
&=&\frac{\|g\|_{\infty}(1-|x|^{2})^{2}}{n-2}\sum_{k=0}^{\infty}
\frac{2}{n}\big(k+\frac{n}{2}\big)(k+1)\frac{(n-2)}{2(k+1)(2k+n)}|x|^{2k}\\
\nonumber
&=&\frac{\|g\|_{\infty}(1-|x|^{2})^{2}}{2n}\sum_{k=0}^{\infty}|x|^{2k}\\
&=&\frac{\|g\|_{\infty}(1-|x|^{2})}{2n}.
\eeq

For $\rho\in(0,1)$, let $F(x)=\mathcal{P}_{f}(\rho x)$, $x\in\mathbb{B}^{n}$. Then
\be\label{eq-ex2}
F(x)=\int_{\mathbb{S}^{n-1}}P(x,\zeta)F(\zeta)d\sigma(\zeta),
\ee
which, together with (\ref{eq-K}), yields that
\be\label{eq-10}
\left|F(x)-\frac{1-|x|^{2}}{(1+|x|^{2})^{\frac{n}{2}}}F(0)\right|
\leq\|\mathcal{P}_{f}\|_{\infty}U(|x|N).
\ee
Applying (\ref{eq-10}), we see that
\beq\label{eq-11}\nonumber
\left|\mathcal{P}_{f}(x)-\frac{1-|x|^{2}}{(1+|x|^{2})^{\frac{n}{2}}}\mathcal{P}_{f}(0)
\right|&=&\lim_{\rho\rightarrow1^{-}}\left|F(x)-\frac{1-|x|^{2}}{(1+|x|^{2})^{\frac{n}{2}}}F(0)\right|\\
&\leq& \|\mathcal{P}_{f}\|_{\infty}U(|x|N).\eeq Hence, by
(\ref{eq-7}) and (\ref{eq-11}), we conclude that
\begin{eqnarray*}
\left|f(x)-\frac{1-|x|^{2}}{(1+|x|^{2})^{\frac{n}{2}}}\mathcal{P}_{f}(0)\right|&=&\left|\mathcal{P}_{f}(x)-
\frac{1-|x|^{2}}{(1+|x|^{2})^{\frac{n}{2}}}\mathcal{P}_{f}(0)-\mathcal{G}_{g}(x)\right|\\
&\leq&\left|\mathcal{P}_{f}(x)-\frac{1-|x|^{2}}{(1+|x|^{2})^{\frac{n}{2}}}\mathcal{P}_{f}(0)\right|+|\mathcal{G}_{g}(x)|\\
&\leq&\|\mathcal{P}_{f}\|_{\infty}U(|x|N)+\frac{\|g\|_{\infty}}{2n}(1-|x|^{2}).
\end{eqnarray*}

Now we prove the sharpness part. For
$x\in\overline{\mathbb{B}^{n}}$, let $g(x)=(-2nM,0,\ldots,0)$ and
$f(x)=(M(1-|x|^{2}),0,\ldots,0)$, where $M$ is a positive constant.
If $x\in\mathbb{S}^{n-1}$, then the optimality of (\ref{eq-thm1}) is
obvious. On the other hand, if $x=0$, then
$$\|\mathcal{P}_{f}\|_{\infty}=0~\mbox{and}~M=|f(0)-\mathcal{P}_{f}(0)|=|f(0)|
= \frac{\|g\|_{\infty}}{2n}=M.
$$
The proof of the theorem is complete.
\qed

\subsection*{Proof of Corollary \ref{cor-1}}
For $\rho\in(0,1)$, let $F(x)=\mathcal{P}_{f}(\rho x)$,
$x\in\mathbb{B}^{n}$. Then, following the proof of Theorem \ref{thm-1} and considering $F$ described by \eqref{eq-ex2},
we deduce that
\beqq
\nonumber
\left|F(x)-\frac{1-|x|}{(1+|x|)^{n-1}}F(0)\right|&=&
\left|\int_{\mathbb{S}^{n-1}}\Big[P(x,\zeta)-\frac{1-|x|}{(1+|x|)^{n-1}}\Big]F(\zeta)d\sigma(\zeta)\right|\\
\nonumber
&\leq&\int_{\mathbb{S}^{n-1}}\Big[P(x,\zeta)-\frac{1-|x|}{(1+|x|)^{n-1}}\Big]|F(\zeta)|d\sigma(\zeta)\\
&\leq&
\|\mathcal{P}_{f}\|_{\infty}\left[1-\frac{1-|x|}{(1+|x|)^{n-1}}\right]
\eeqq
and therefore,
\beqq
\nonumber
\left|\mathcal{P}_{f}(x)-\frac{1-|x|}{(1+|x|)^{n-1}}\mathcal{P}_{f}(0)
\right|&=&\lim_{\rho\rightarrow1^{-}}\left|F(x)-\frac{1-|x|}{(1+|x|)^{n-1}}F(0)\right|\\
&\leq& \|\mathcal{P}_{f}\|_{\infty}\left[1-\frac{1-|x|}{(1+|x|)^{n-1}}\right]
\eeqq
which by  Theorem \ref{thm-1} leads to
\begin{eqnarray*}
\left|f(x)-\frac{1-|x|}{(1+|x|)^{n-1}}\mathcal{P}_{f}(0)\right|&=&\left|\mathcal{P}_{f}(x)-\frac{1-|x|}{(1+|x|)^{n-1}}\mathcal{P}_{f}(0)-\mathcal{G}_{g}(x)\right|\\
&\leq&\left|\mathcal{P}_{f}(x)-\frac{1-|x|}{(1+|x|)^{n-1}}\mathcal{P}_{f}(0)\right|+|\mathcal{G}_{g}(x)|\\
&\leq&\|\mathcal{P}_{f}\|_{\infty}\left[1-\frac{1-|x|}{(1+|x|)^{n-1}}\right]
+\frac{\|g\|_{\infty}}{2n}(1-|x|^{2}).
\end{eqnarray*}
The proof of the corollary is complete.
 \qed

\begin{Lem}{\rm (\cite[Lemma 2.3]{K5})}\label{Lem-A}
For $r\in[0,1]$, let $\varphi(r)=\frac{\partial U(rN)}{\partial r}$.
Then $\varphi(r)$ is decreasing on $r\in[0,1]$ and
$$\varphi(r)\geq \left . \frac{\partial U(rN)}{\partial r} \right |_{r=1}=\varphi(1)=A_{n},
$$
where $A_{n}$ is the same as in
Theorem {\rm\ref{thm-2}}.
\end{Lem}

\subsection*{Proof of Theorem \ref{thm-2}}
For $x\in\mathbb{B}^{n}$, there is a $\rho\in(r,1)$ such that
\be\label{eq-12c} \frac{1-U(rN)}{1-|x|}=\frac{\partial U(\rho
N)}{\partial r}, \ee where $r=|x|$. Now, for $n\geq3$ and a given
$g\in\mathcal{C}(\overline{\mathbb{B}^{n}},\mathbb{R}^{m})$, by
(\ref{eq-p}) with $f$ in place of $\psi$, we have
$$f(x)=\mathcal{P}_{f}(x)-\mathcal{G}_{g}(x),
$$
where
$\mathcal{P}_{f}$ and $\mathcal{G}_{g}$ are defined as in \eqref{eq-p1a}. Since
$f(0)=\mathcal{P}_{f}(0)-\mathcal{G}_{g}(0)=0$, by (\ref{eq-7}),
Theorem \ref{thm-1} and the assumptions, we see that
\beq\label{eq-13c}\nonumber
|f(\zeta)-f(r\zeta)|&=&\left|f(\zeta)+\mathcal{P}_{f}(0)\frac{1-r^{2}}{(1+r^{2})^{\frac{n}{2}}}-
\mathcal{G}_{g}(0)\frac{1-r^{2}}{(1+r^{2})^{\frac{n}{2}}}-f(r\zeta)\right|\\
\nonumber &\geq&1-
\left|f(r\zeta)-\mathcal{P}_{f}(0)\frac{1-r^{2}}{(1+r^{2})^{\frac{n}{2}}}\right|-|\mathcal{G}_{g}(0)|\frac{1-r^{2}}{(1+r^{2})^{\frac{n}{2}}}\\
&\geq&1-U(rN)-\frac{\|g\|_{\infty}}{2n}(1-r^{2})-\frac{\|g\|_{\infty}}{2n}\frac{1-r^{2}}{(1+r^{2})^{\frac{n}{2}}},
\eeq
where $r\in[0,1)$. Finally, by (\ref{eq-12c}), (\ref{eq-13c}) and Lemma \Ref{Lem-A}, there is a $\rho\in(r,1)$ such that
\begin{eqnarray*}
\frac{|f(\zeta)-f(r\zeta)|}{1-r}
&\geq&\frac{1-U(rN)}{1-r}-\frac{\|g\|_{\infty}}{2n}(1+r)-\frac{\|g\|_{\infty}}{2n}\frac{1+r}{(1+r^{2})^{\frac{n}{2}}}\\
&\geq&\frac{\partial U(\rho N)}{\partial r}-\frac{\|g\|_{\infty}}{2n}(1+r)\left (1+\frac{1}{(1+r^{2})^{\frac{n}{2}}}\right )\\
&\geq&A_{n}-\frac{\|g\|_{\infty}}{2n}(1+r)\left (1+\frac{1}{(1+r^{2})^{\frac{n}{2}}}\right ),
\end{eqnarray*}
which gives \eqref{eq-Sch}.
The sharpness part easily follows from \cite[Theorem 2.5]{K5}. The
proof of the theorem is complete. \qed

\begin{Thm}{\rm (\cite[Theorem 2.12]{MM})}\label{Lem-MM}
Let $u$ be a bounded harmonic function from $\mathbb{B}^{n}$ into
$\mathbb{R}$, where $n\geq3$. Then, for $x\in\mathbb{B}^{n}$,
$$|\nabla
u(x)|\leq\frac{\|u\|_{\infty}}{1-|x|}\sup_{\gamma>0}C(\gamma,|x|),$$
where $C(\gamma,|x|)$ is defined in Theorem {\rm\ref{thm-3}}.
\end{Thm}

\subsection*{Proof of Theorem \ref{thm-3}}
Let $n\geq3$ and $g\in\mathcal{C}(\overline{\mathbb{B}^{n}},\mathbb{R}^{m})$ be given. As before, by
(\ref{eq-p}), we have
$$f(x)=\mathcal{P}_{f}(x)-\mathcal{G}_{g}(x),
$$
where $\mathcal{P}_{f}$ and $\mathcal{G}_{g}$ are defined as in
\eqref{eq-p1a}. 
Set $\mathcal{G}_{g}=(\mathcal{G}_{g,1},\ldots,\mathcal{G}_{g,m})$
and $g=(g_{1},\ldots,g_{m})$. For $k\in\{1,2,\ldots,m\}$, we let
$$I_{k}=\int_{\mathbb{B}^{n}}\left|\frac{x-y}{|x-y|^{n}}-\frac{|y|^{2}x-y}{[x,y]^{n}}\right||g_{k}(y)|dV(y).
$$
If we apply Cauchy-Schwarz's inequality and \cite[Theorem 2.1]{K2}, it follows that
\beq\label{eq16c}
I_{k}^{2}&\leq&
\int_{\mathbb{B}^{n}}\left|\frac{x-y}{|x-y|^{n}}-\frac{|y|^{2}x-y}{[x,y]^{n}}\right|dV(y)\\
\nonumber &&\times\bigg[\int_{\mathbb{B}^{n}}
\left|\frac{x-y}{|x-y|^{n}}-\frac{|y|^{2}x-y}{[x,y]^{n}}\right||g_{k}(y)|^{2}dV(y)\bigg]\\
\nonumber &\leq&
\frac{2n\pi^{\frac{n}{2}}}{(n+1)\Gamma\big(\frac{n}{2}\big)}\int_{\mathbb{B}^{n}}
\left|\frac{x-y}{|x-y|^{n}}-\frac{|y|^{2}x-y}{[x,y]^{n}}\right||g_{k}(y)|^{2}dV(y),\eeq
which yields that
\beq\label{eq17c}|\nabla
\mathcal{G}_{g,k}(x)|^{2}&\leq&\frac{I_{k}^{2}}{\omega_{n-1}^{2}}\nonumber \\
&\leq&\frac{n}{(n+1)\omega_{n-1}}\int_{\mathbb{B}^{n}}
\left|\frac{x-y}{|x-y|^{n}}-\frac{|y|^{2}x-y}{[x,y]^{n}}\right||g_{k}(y)|^{2}dV(y).
\eeq
Then, using (\ref{eq17c}) and \cite[Theorem 2.1]{K2}, we obtain
\beq\label{eq-15c}\nonumber
|D_{\mathcal{G}_{g}}(x)|&=&\sup_{\theta\in\mathbb{S}^{n-1}}\left(\sum_{k=1}^{m}|\langle\nabla
\mathcal{G}_{g,k}(x),\theta\rangle|^{2}\right)^{\frac{1}{2}}\\
\nonumber &\leq&\left(\sum_{k=1}^{m}|\nabla
\mathcal{G}_{g,k}(x)|^{2}\right)^{\frac{1}{2}}\\ \nonumber
&\leq&\bigg[\frac{n}{(n+1)\omega_{n-1}}\int_{\mathbb{B}^{n}}
\left|\frac{x-y}{|x-y|^{n}}-\frac{|y|^{2}x-y}{[x,y]^{n}}\right|\sum_{k=1}^{m}|g_{k}(y)|^{2}dV(y)\bigg]^{\frac{1}{2}}\\
\nonumber&\leq&\left[\frac{n}{(n+1)\omega_{n-1}}\right]^{\frac{1}{2}}\|g\|_{\infty}\left(\int_{\mathbb{B}^{n}}
\left|\frac{x-y}{|x-y|^{n}}-\frac{|y|^{2}x-y}{[x,y]^{n}}\right|dV(y)\right)^{\frac{1}{2}}\\
\nonumber
&\leq&\left[\frac{n}{(n+1)\omega_{n-1}}\right]^{\frac{1}{2}}\left(\frac{n\omega_{n-1}}{n+1}\right)^{\frac{1}{2}}\|g\|_{\infty}\\
&=&\frac{n}{n+1}\|g\|_{\infty}.
\eeq

Now we estimate $|D_{\mathcal{P}_{f}}|$. For $x\in\mathbb{B}^{n}$, we may let
$\mathcal{P}_{f}(x)=(\mathcal{P}_{f,1}(x),\ldots,\mathcal{P}_{f,m}(x)).$ Then, for any
$\theta\in\mathbb{S}^{n-1}$ and $k\in\{1,2,\ldots,m\}$, by
Cauchy-Schwarz's inequality, we have
\begin{eqnarray*}
\left|\langle\nabla
\mathcal{P}_{f,k}(x),\theta\rangle\right|^{2}&=&\left|\int_{\mathbb{S}^{n-1}}\langle\nabla
P(x,\zeta),\theta\rangle\mathcal{P}_{f,k}(\zeta)d\sigma(\zeta)\right|^{2}\\
 &\leq&\left(\int_{\mathbb{S}^{n-1}}\big|\langle\nabla
P(x,\zeta),\theta\rangle\big||\mathcal{P}_{f,k}(\zeta)|d\sigma(\zeta)\right)^{2}\\
&\leq&\delta(x,\theta)\int_{\mathbb{S}^{n-1}}\big|\langle\nabla
P(x,\zeta),\theta\rangle\big||\mathcal{P}_{f,k}(\zeta)|^{2}d\sigma(\zeta),
\end{eqnarray*}
which gives that
\beq\label{eq-cp1}\nonumber
 \left(\sum_{k=1}^{m}\big|\langle\nabla
\mathcal{P}_{f,k}(x),\theta\rangle\big|^{2}\right)^{\frac{1}{2}}&\leq&(\delta(x,\theta))^{\frac{1}{2}}
\left(\int_{\mathbb{S}^{n-1}}\big|\langle\nabla
P(x,\zeta),\theta\rangle\big|\sum_{k=1}^{m}|\mathcal{P}_{f,k}(\zeta)|^{2}d\sigma(\zeta)\right)^{\frac{1}{2}}\\
&\leq&\delta(x,\theta)\|\mathcal{P}_{f}\|_{\infty},
\eeq
where
$$\delta(x,\theta)=\int_{\mathbb{S}^{n-1}}\big|\langle\nabla P(x,\zeta),\theta\rangle\big|d\sigma(\zeta).
$$
Applying (\ref{eq-cp1}), \cite[Lemma 2.3]{MM} and Theorem \Ref{Lem-MM}, we
see that, for $x\in\mathbb{B}^{n}$,
\beq\label{eq-cp2}
\nonumber|D_{\mathcal{P}_{f}}(x)|&=&\sup_{\theta\in\mathbb{S}^{n-1}}\left(\sum_{k=1}^{m}\big|\langle\nabla
\mathcal{P}_{f,k}(x),\theta\rangle\big|^{2}\right)^{\frac{1}{2}}\\
\nonumber &\leq&\|\mathcal{P}_{f}\|_{\infty}
\sup_{\theta\in\mathbb{S}^{n-1}}\delta(x,\theta)\\
&\leq&
\frac{\|\mathcal{P}_{f}\|_{\infty}}{1-|x|}\sup_{\gamma>0}C(\gamma,|x|)
\eeq
By   (\ref{eq-15c}) and (\ref{eq-cp2}),   we conclude that
$$|D_{f}(x)|\leq|D_{\mathcal{P}_{f}}(x)|+|D_{\mathcal{G}_{g}}(x)|
\leq\frac{\|\mathcal{P}_{f}\|_{\infty}}{1-|x|}\sup_{\gamma>0}C(\gamma,|x|)+ \frac{n}{n+1}\|g\|_{\infty}, ~~x\in\mathbb{B}^{n}.
$$
The proof the theorem is complete. \qed


\section{An application of the Schwarz Lemma}\label{csw-sec3}

Let $f:\,\overline{\Omega}\to
\mathbb{R}^n$ be a differentiable mapping and $x$ be a regular value
of $f$, where $x\notin f(\partial\Omega)$ and $\Omega\subset
\mathbb{R}^n$ is a bounded domain. Then the degree
$\deg(f,\Omega,x)$ is defined by the formula (cf. \cite{Ll,V})
$$\deg(f,\Omega,x):=\sum_{y\in f^{-1}(x)\cap\Omega}\mbox{ sign} \big(\det J_{f} (y)\big).
$$

\begin{Lem}\label{Lem-X}
The $\deg(f,\Omega,x)$ satisfies the following properties {\rm (cf.
\cite[p.~125-129]{RR}):}
\begin{enumerate}
\item[(a)] If $x\in\mathbb{R}^ n \backslash f(\partial D)$ and $\deg(f,\overline{\Omega},x)\neq 0$, then there exists an $w\in\Omega$ such that $f(w)=x$.
\item[(b)]\label{(II)} If $D$ is a domain with $\overline{D}\subset\Omega$ and $x\in \mathbb{R}^ n \backslash f(\partial D)$,
then $\deg(f,D,x)$ is a constant on each component of $\mathbb{R}^
n \backslash f(\partial D)$.
\end{enumerate}
\end{Lem}

\begin{lem}\label{lem-L}
Let $n\geq3$ and $m\geq1$. For a given
$g\in\mathcal{C}(\overline{\mathbb{B}^{n}},\mathbb{R}^{m})$, if
$f\in\mathcal{C}^{2}(\mathbb{B}^{n},\mathbb{R}^{m})\cap\mathcal{C}(\mathbb{S}^{n-1},\mathbb{R}^{m})$
satisfies $\Delta f=g$ and $\|f\|_{\infty}+\|g\|_{\infty}\leq M$ for some constant $M>0$,
then there is a constant $L>0$ such that
$$|f(x_{1})-f(x_{2})|\leq L|x_{1}-x_{2}| ~\mbox{ for all
$x_{1},x_{2}\in\overline{\mathbb{B}^{n}(x_{0},\rho_{0})}$,}
$$
where $x_{0}\in\mathbb{B}^{n}$ and
$\rho_{0}\in(0,1-|x_{0}|)$ are some constants.
\end{lem}
\bpf By Theorem \ref{thm-3} or (\ref{eq-19cp}), for all
$x\in\overline{\mathbb{B}^{n}(x_{0},\rho_{0})}$, we have
$$|D_{f}(x)|\leq\frac{\|\mathcal{P}_{f}\|_{\infty}}{1-|x|}\sup_{\gamma>0}C(\gamma,|x|)+
\frac{n}{n+1}\|g\|_{\infty}\leq\frac{2nM}{1-\rho_{0}-|x_{0}|}+\frac{n}{n+1}M:=L,
$$
which implies, for all $x_{1},x_{2}\in\overline{\mathbb{B}^{n}(x_{0},\rho_{0})}$,
$$|f(x_{1})-f(x_{2})|\leq\int_{[x_{1},x_{2}]}|D_{f}(x)|\, |dx|\leq L\int_{[x_{1},x_{2}]}|dx|=L|x_{1}-x_{2}|,$$
where $[x_{1},x_{2}]$ is the segment from $x_{1}$ to $x_{2}$ (or
$x_{2}$ to $x_{1}$) with the endpoints $x_{1}$ and $x_{2}$. \epf

\subsection*{Proof of Theorem \ref{thm-4}}

We prove the theorem by the method of contradiction. Suppose that
the result is not true. Then there is a sequence $\{a_k\}$ of
positive real numbers, and a sequence of functions $\{f_k\}$ with
$f_{k}\in\mathcal{F}_{g}^{M}$, such that $\{a_k\}$ tends to $0$ and
$a_k \notin f_k(\mathbb{B}^{n})$ for $k\in\{1,2,\ldots\}$. By
Theorem \ref{thm-1}, Lemma \ref{lem-L} and Arzel${\rm \grave{a}}$-Ascoli's theorem, 
we know that there is a subsequence $\{f_{k}^{\ast}\}$ of $\{f_{k}\}$ which converges
uniformly on a compact subset of $\mathbb{B}^{n}$ to a function
$f^{\ast} $. For each $k$, it is easy to see that the function
$h_k=f_{k}^{\ast}-f_{1}^{\ast}$ is harmonic. As a consequence, the
sequence $\{h_k\}$ converges uniformly on compact subsets of
$\mathbb{B}^{n}$ to $f^{\ast}-f_{1}^{\ast}$ and therefore, the
partial derivatives of $f_{k}^{\ast}$ converge uniformly on compact
subsets of $\mathbb{B}^{n}$ to the partial derivatives of $
f^{\ast}$. In particular, $f_{k}^{\ast} (0) \rightarrow f^{\ast}(0)$
and $J_{f_{k}^{\ast}} (0) \rightarrow J_{f^{\ast}} (0)$ as
$k\rightarrow\infty$, which imply that $f^{\ast} \in
\mathcal{F}_{g}^{M}$. Since $J_{f^{\ast}} (0)-1=|f^{\ast}(0)|=0,$
there are $r_0\in(0,1)$ and $c_1
> 0$ such that $J_{f^{\ast}} >0$ on $
\overline{\mathbb{B}^{n}(0,r_{0})}$,
$f^{\ast}(\mathbb{B}^{n}(0,r_{0})) \supset
\overline{\mathbb{B}^{n}(0,c_1)}$ and $|f^{\ast}(x)| \geq c_1$ for
$x \in
\partial\mathbb{B}^{n}(0,r_{0})$.

Now, we let $c_2=c_1/2$, $\mathbb{B}_{r_{0}}=\mathbb{B}^{n}(0,r_{0})$ and
$\mathbb{B}_{c_{2}}= \mathbb{B}^{n}(0,c_2) $. Then there is a $k_0$
such that $|f^{\ast}_{k}(x)| \geq c_2$ for $k \geq k_0$ and
$J_{f^{\ast}_k}
>0$ on $\overline{\mathbb{B}_{r_{0}}}$. Since $\deg(f^{\ast}_k,\mathbb{B}_{r_{0}},0)\geq 1$,
by Lemma \Ref{Lem-X}, we see that, for $y \in \mathbb{B}_{c_{2}}$
and $k \geq k_0$, $\deg(f^{\ast}_k,\mathbb{B}_{r_{0}},y)\geq 1$.
Hence, for $k \geq k_0$, $f^{\ast}_k(\mathbb{B}_{r_{0}})\supset
\mathbb{B}_{c_{2}}$ which contradicts our assumption.
The proof of the theorem is complete.
\qed

\bigskip

{\bf Acknowledgements:}  This research was partly supported by the
National Natural Science Foundation of China ( No. 11571216 and No.
11401184),  the Science and Technology Plan Project of Hunan
Province (No. 2016TP1020) and the Construct Program of the Key
Discipline in Hunan Province. The second author is currently on leave from IIT Madras.

\subsection*{Conflict of Interests}
The authors declare that there is no conflict of interests regarding the publication of this paper.

\end{document}